\documentclass[10pt]{article}
\usepackage{amsmath,amsfonts,latexsym,amssymb}
\usepackage[colorlinks=true,
pdfstartview=FitV, linkcolor=blue, citecolor=blue,
urlcolor=blue]{hyperref}
\usepackage[latin1]{inputenc}
\newtheorem{theorem}{Theorem}[subsection]
\newtheorem{lemma}[theorem]{Lemma}
\newtheorem{proposition}[theorem]{Proposition}
\newtheorem{corollary}[theorem]{Corollary}

\newcommand{\proof}{\noindent{\sc Proof : }}
\newcommand{\qed}{{\sc Q.e.d.}}
\newcommand{\rmks}{{\sc Remarks:}}
\newcommand{\rmk}{{\sc Remark:}}

\begin{document}

\title{Well displacing representations\\ and orbit maps}
\author{Thomas DELZANT\thanks{IRMA et Université L. Pasteur, 7 rue Descartes, F-67084 Strasbourg}\and Olivier GUICHARD\thanks{CNRS, Orsay cedex, F-91405; Univ. Paris-Sud, Laboratoire de Mathématiques, Orsay F-91405 Cedex; } \and
François LABOURIE \thanks{Univ. Paris-Sud, Laboratoire de Mathématiques, Orsay F-91405 Cedex; CNRS, Orsay cedex, F-91405} \and Shahar MOZES\thanks{Institute of Mathematics, Hebrew University, Jerusalem, Israel}}
\maketitle
\section{Introduction}

\subsection{Displacement function}\label{defwd}
Let $\gamma$ be an isometry of a metric space $X$. We recall that the {\em displacement} of $\gamma$ is
$$
d_X(\gamma)=\inf_{x\in X}d(x,\gamma(x)).
$$
We observe that $d_{X}(\gamma)$ is an invariant of the conjugacy class of $\gamma$.

If $C$ is the Cayley graph of a group $\Gamma$ with set of generators $S$ and word length $\Vert\ \Vert_S$, the displacement function is called the {\em translation length} and is denoted by $\ell$
$$
\ell_S(\gamma)=\inf_{\eta}\Vert\eta\gamma\eta^{-1}\Vert_S.
$$
We finally say the action by isometries on $X$ of a group  $\Gamma$ is {\em well displacing}, if given a set $S$ of generators of $\Gamma$, there exist positive constants $A$ and $B$ such that
$$
d_X(\gamma)\geqslant  A \ell_S(\gamma)-B.
$$
This definition does not depends on the choice of $S$. As first examples,
it is easy to check that cocompact groups  are well displacing, as well as convex cocompact whenever $X$ is {\em Hadamard} (i.e. complete, non positively curved and simply connected). We recall that a cocompact action is by definition a properly discontinuous action whose quotient is compact, and a convex cocompact is an action such that there exists a convex invariant on which the action is cocompact.

The notion naturally arose in  \cite{Labourie:2005a} where it is shown  that for well displacing representations, the energy functional is proper on Teichmüller space and that moreover a large class of representations of surface groups are well displacing.

The purpose  of this note is to investigate the relation of well displacing actions, to actions for which  the orbit maps are quasi isometric embeddings. For hyperbolic groups the two notions turn out to be equivalent. However in general, the answer is a little more involved.
\begin{itemize}
\item There exist infinite groups whose actions are always well displacing. However, an infinite group do have actions for which orbit maps are not quasi isometric embeddings: the trivial action. {\em Therefore, in general, well displacing representation do not have orbit maps which are quasi isometric embeddings.}
\item{\em There exist a class of group for which if an orbit map is a quasi isometric embedding then the action is well displacing}. The main example is the class of hyperbolic groups. In general these are the groups, for which the stable length is comparable to to the translation length ({\em cf.} Proposition \ref{prop:qi->wd}).
\item {\em There exists a class of group for which every well displacing action have orbit maps which are quasi isometric embeddings}. This class of group are those which are {\em undistorted in their conjugacy classes} (see  Section    \ref{def:U}) -- or $U$-property -- reminiscent of a statement of Abels, Margulis and Soifer about proximal maps \cite{Abels:1995}. We prove in Theorem \ref{thm_UCC}, that some class of linear groups  -- in particular lattices --- enjoys the $U$-property. We also prove that hyperbolic groups have the $U$-property.
\item {\em There exists action of groups -- which have the $U$-property -- whose orbit maps are quasi isometric embeddings but which are not well displacing}. The  simplest example is ${\rm SL}(n,\mathbb Z)$ acting on ${\rm SL}(n,\mathbb R)/SO(n,\mathbb R)$. More generally in Corollary \ref{cayl2}, we show that  if a group  has {\em infinite contortion} (see Section \ref{def:contort}) -- in particular every residually finite group -- then any linear  representation which contains a unipotent is not well displacing.
\end{itemize}

Finally, we investigate in a complement a property related to infinite contortion which we call the bounded depth roots property (see Section  \ref{bdr}). We finish by a simple question :
 Do all linear groups have $U$-property?

 \tableofcontents

\section{Infinite groups whose actions are always well displacing}
We have
\begin{proposition}
There exists infinite finitely generated groups whose actions are always well displacing. Hence there exists action for which the orbit maps are not quasi isometric embeddings, but which are well displacing.
\end{proposition}
\proof Denis Osin  \cite{Osin:2004} has constructed infinite finitely generated groups with exactly $n$ conjugacy classes. Any action of such a group is well displacing. For the second part, we just take the trivial action  on a point.
\section{Orbit maps which are quasi isometry implies well displacing}
If a group $\Gamma$ acts by isometries on a metric space $X$, we define the {\em stable norm with respect to $X$} by
$$
[g]^X_\infty=\liminf_{n\rightarrow\infty}\frac{1}{n}d(x_0,g^n(x_0)).
$$
We observe that this quantity does not depends on the choice of the base point $x$. The {\em stable norm} $[g]_\infty$ is the stable norm with respect to the Cayley graph of $\Gamma$.

We now prove the following easy result

\begin{proposition}\label{prop:qi->wd}
If a group  $\Gamma$
is such that there exists $\alpha >0$ such that
$$
\forall g\in\Gamma, \ \  [g]_\infty\geqslant \alpha. \ell (g),
$$
then every action of $\Gamma$ on $(X,d)$ for which the orbit map is a quasi isometric embedding is well displacing.
\end{proposition}

\proof by definition, if the orbit map is a quasi isometric embedding, for every $x\in X$ there exists some constant $A$ and $B$ such that
$$
A\Vert \gamma\Vert +B\geqslant d(x,\gamma(x))\geqslant A^{-1}\Vert \gamma\Vert-B.
$$
It follows that
$$
A[\gamma]_\infty\geqslant [\gamma]_\infty^X \geqslant A^{-1}[\gamma]_\infty.
$$
Now we remark that
$$
[\gamma]^X_\infty\leqslant d_X(\gamma).
$$
The result follows.\qed
\vskip 0.3 truecm

\rmks
\begin{itemize}
\item We will show later on that this inequality fails for $SL(n,\mathbb Z)$ for $n\geqslant 3$.
\item On the other hand, we observe following
\cite[p 119]{Coornaert:1990},  that for $\Gamma$ a hyperbolic group, the stable norm of an element coincides up to a constant with its translation length : there exists a contant $K$ such that $\vert \ell(g)-[g]_\infty\vert\leqslant K$.

\end{itemize}
Therefore, we have

\begin{corollary}\label{corohyper}
Let $\Gamma$ be a hyperbolic group. If an isometric action is such that the orbit maps are quasiisometries, then this action is well displacing.
\end{corollary}

\section{Groups whose well displacing actions have orbit maps which are quasi isometric embeddings}
\label{def:U}
We say a finitely generated group is {\em undistorted in its conjugacy classes} -- in short,  has {\em $U$-property} -- , if there exists finitely many elements $g_1,\ldots g_p$ of $\Gamma$, positive constants $A$ and $B$  such that
$$
\forall \gamma\in \Gamma, \Vert \gamma\Vert\leqslant A\sup_{1\leqslant i\leqslant n}\ell (g_i\gamma)+B.
$$
\rmks
\begin{itemize}
\item This property is clearly independent of $S$.
\item This property is satisfied by free groups and commutative groups. On the other hands the groups constructed by D. Osin \cite{Osin:2004} described in the above paragraph do not have this $U$-property.
\item Note also that this
property is very similar to the statement of Abels, Margulis and
Soifer's result \cite[Theorem 4.1]{Abels:1995} and it is no
surprise that their result plays a role in the proof of
Theorem~\ref{thm_UCC}.
\item Finally, by the conjugacy invariance of the translation length, $\ell( g_i
\gamma) = \ell( \gamma g_{i})$ so we will indifferently write
this property with left or right multiplication by the finite family
$(g_i)$.
\end{itemize}

\begin{lemma}
If $\Gamma$ has $U$-property, then every well displacing action have orbit maps which are quasi isometric embeddings.
\end{lemma}
\proof Indeed, assume that $\Gamma$ acts on $X$ by isometries and that the action is well displacing.
In particular we have positive constants $\alpha$ and $\beta$ so that $$
\forall x\in X, \ \ \gamma\in\Gamma, \ \ d(x,\gamma(x))\geqslant  \alpha \ell_S(\gamma)-\beta.
$$
Moreover
there exists finitely many elements $g_1,\ldots g_p$ of $\Gamma$, positive constants $A$ and $B$  such that
$$
\forall \gamma\in \Gamma, A\Vert \gamma\Vert-B\leqslant \sup_{1\leqslant i\leqslant n}\ell (g_i\gamma).
$$
Hence, let $x\in X$, then
\begin{eqnarray*}
d(x,\gamma(x))&\geqslant & \sup_{1\leqslant i\leqslant n} d(x,g_i.\gamma (x))-\sup_{1\leqslant i\leqslant n}d(x,g_i (x))
\\
&\geqslant & \alpha \sup_{1\leqslant i\leqslant n}\ell(g_i\gamma)-
\beta- \sup_{1\leqslant i\leqslant n}d(x,g_i (x))\\
&\geqslant & \alpha.{A} \Vert \gamma\Vert -B\alpha -
\beta - \sup_{1\leqslant i\leqslant n}d(x,g_i(x)).
\end{eqnarray*}
Hence the orbit map is a quasi isometric embedding.
\qed

We will prove in the next section,
\begin{theorem}  \label{thm_linG}
Every uniform lattice -- and non uniform lattice in higher rank -- in characteristic zero has  $U$-property. In particular every surface group has $U$-property.
\end{theorem}
Moreover, we show that hyperbolic groups have $U$-property.

\subsection{Linear Groups having $U$-property}
We prove the following result which implies Theorem \ref{thm_linG}

\begin{theorem}
  \label{thm_UCC}
  Let $\Gamma$ be a  finitely  generated group and $\mathbf{G}$ a reductive
  group defined over a field $F$, suppose that
  \begin{itemize}
  \item there exists a homomorphism $\Gamma \rightarrow \mathbf{G}(F)$ with Zariski dense image,
  \item there are a finitely many  field homomorphisms $(i_v)_{v\in S}$  of $F$ in local fields $F_v$  such that the diagonal embedding $\Gamma \rightarrow \Pi_{v \in
      S} \mathbf{G}( F_v)$ is a quasi-isometric embedding.
  \end{itemize}
  Then the group $\Gamma$ has $U$-property.
\end{theorem}

Since lattices are Zariski dense (Borel theorem) and that higher rank
irreducible lattices
in characteristic zero are quasi isometrically embedded
(\cite{Lubotzky:2000}), this theorem implies Theorem \ref{thm_linG} for higher
rank lattices. The same holds for all uniform lattices.
A specific corollary is the following

\begin{corollary}
Let $\Gamma$ be a finitely generated that is quasi isometrically embedded and Zariski dense in a reductive group $G(F)$ where $F$ is a local field. Then $\Gamma$ has $U$-property.
\end{corollary}

\subsubsection{Generalities for $U$-property}
\label{sec_gen}

We prove two  lemmas for the $U$-property.

\begin{lemma}
  \label{lem_normalsubgroup}. Let $\Gamma$ be a finitely generated group. Let $\Gamma_{0} \vartriangleleft \Gamma$ be  a  normal subgroup of finite
  index. If $\Gamma_{0}$ has  $U$-property, so  has  $\Gamma$.
\end{lemma}

We do not know whether the converse statement holds, in other words whether $U$-property is a property of commensurability classes.
\vskip 0.2cm

\proof We first observe that every finite index subgroup of a finitely generated group is finitely generated.
Let $S_{0}$ be a generating set for $\Gamma_{0}$ and write $\Gamma$ as the
union of left cosets for $\Gamma_{0}$
\begin{equation*}
  \Gamma = \bigcup_{t \in T} \Gamma_{0} \cdot t.
\end{equation*}
We assume that $T$ is symmetric.
Clearly $S= S_{0} \cup T$ is a generating set for $\Gamma$.

We denote
${\Vert}\cdot{\Vert}_{ \Gamma_{0}}$ and $\ell_{\Gamma_0}$ the word  and translation lengths for $\Gamma_{0}$.

We observe that $\Gamma_0$ is quasi-isometrically embedded in $\Gamma$. Hence there exist positive  constants $\alpha$ and $\beta$ such that
\begin{equation}
\forall\gamma\in\Gamma_0, \ \ \Vert \gamma\Vert _{\Gamma_0}\geqslant  \Vert \gamma\Vert _{\Gamma}\geqslant  \alpha\Vert \gamma\Vert _{\Gamma_0}-\beta.\label{ineq:quasi}
\end{equation}

For any $\gamma$ in $\Gamma$, we write $\gamma  = \gamma_{0} t_0$ with
$t_0 \in T$ and $\gamma_{0} \in \Gamma_{0}$. Hence
\begin{equation*}
  {\Vert} \gamma{\Vert}_{ \Gamma} \leqslant {\Vert} \gamma_{0}{\Vert}_{ \Gamma_{0}} + 1 \leqslant A \sup
  \ell_{ \Gamma_{0}}( \gamma_{0} g_i) +B+1 \label{ineq:22}
\end{equation*}
since $\Gamma_{0}$ has $U$-property.

Finally we need to compare $\ell_{ \Gamma}$ and  $\ell_{ \Gamma_0}$. Let $\delta$ in $\Gamma_{0}$, then
\begin{eqnarray}
  \ell_{ \Gamma}( \delta) &=& \inf_{ t \in T, \eta \in \Gamma_0} {\Vert} t \eta t \delta \eta^{-1} t^{-1}
  {\Vert}_{\Gamma}  \cr
  &\geqslant & \inf_{\eta \in \Gamma_0} {\Vert} \eta \delta \eta^{-1}\Vert_{\Gamma} -2\cr
 &\geqslant & \alpha\inf_{\eta \in
    \Gamma_{0}} {\Vert} \eta \delta \eta^{-1} {\Vert}_{\Gamma_{0}}-\beta-2\cr
    & \geqslant &
\alpha  \ell_{ \Gamma_{0}}( \delta )-\beta -2\label{ineq:21}
\end{eqnarray}

Finally combining Inequalities (\ref{ineq:22}) and (\ref{ineq:21}) , we have
\begin{equation*}
  {\Vert} \gamma{\Vert}_{\Gamma} \leqslant \frac{A}{\alpha}\sup_{t\in T, i
    \in I} \ell_{\Gamma}( \gamma t^{-1} \beta_i) + B+1  + \frac{\beta+2}{\alpha}.
\end{equation*}
This is exactly the $U$-property for $\Gamma$.
\qed

Also
\begin{lemma}
  \label{lem_finiteindex} Let  $\Gamma$ be a  finitely generated group.
  Suppose that $\Gamma \rightarrow \Gamma_{0}$ is onto with finite
  kernel.  Then the group $\Gamma$ has $U$-property if and only if $\Gamma_{0}$ has.
\end{lemma}

\proof We choose a generating set $S$ for $\Gamma$ which contains the kernel
of $\Gamma \rightarrow \Gamma_{0}$. We choose the generating set $S_0$ for $\Gamma_{0}$
to be the image of $S$. Then we have, using surjectivity, for all
$\gamma$ projecting to $\gamma_{0}$
\begin{eqnarray*}
  {\Vert} \gamma{\Vert}_\Gamma \geqslant  & {\Vert} \gamma_{0}{\Vert}_{ \Gamma_{0}}& \geqslant  {\Vert}\gamma{\Vert}_{\Gamma}
  -1,\\
 \ell_\Gamma( \gamma) \geqslant  &\ell_{\Gamma_{0}}( \gamma_{0})& \geqslant
  \ell_{\Gamma}( \gamma) -1.
\end{eqnarray*}
These two inequalities  enable us to transfer the $U$-property  from $\Gamma$ to $\Gamma_{0}$ and vice versa.
\qed

\subsubsection{Proximality}
\label{sec_prox}

We recall the notion of proximality and a result of Abels, Margulis
and Soifer.

Let $k$ be a local field.  Let $V$ be a finite dimensional $k$-vector space equipped with a norm. Let $d$ be the induced metric on $\mathbb{P}(V)$. Let $r$ and $\epsilon$ be positive numbers such that
$$r>2\epsilon.$$

An element $g$ of $\textup{{\rm SL}}(V)$ is
said to be $(r, \varepsilon)$-\emph{proximal},  if there exist  a point $x_+$ in
$\mathbb{P}(V)$ and an hyperplane $H$ in $V$  such that
\begin{itemize}
\item $d(x_+,\mathbb{P}( H))\geqslant  r$,
\item $\forall x\in\mathbb{P}(V), \ \  d( x, \mathbb{P}(H) ) \geqslant  \varepsilon\implies d ( g\cdot x, x_+ ) \leqslant \varepsilon$.
\end{itemize}
In particular a proximal element has a unique eigenvalue of highest
norm. Conversely, if  an element $g$ admits a unique eigenvalue of
highest norm, then  some power of $g$ is proximal (for some $(r,\epsilon)$).

We cite the needed result from \cite{Abels:1995} and
\cite{Benoist:1997}.

\begin{theorem}[\cite{Abels:1995} Theorem~5.17]
  \label{coro_AMS}
Let $\mathbf{G}$ a semisimple group over a field $F$. Let $(i_v)_{v\in V}$ be  finitely many  field homomorphisms  of $F$ in local fields $F_v$. Let $\rho_v: \mathbf{G}( F_v) \rightarrow \mathrm{GL}( n_v , F_v)$ be an
  irreducible representation of $\mathbf{G}(F_v)$ for each $v$.

  Suppose that $\Gamma$ is a Zariski dense subgroup of $\mathbf{G}(
  F)$.  Suppose that for every $v$, $\rho_v(\Gamma)$ contains  proximal elements. Then there exist
 \begin{itemize}
 \item $r>2 \varepsilon>0$
 \item  a
  finite subset $\Delta \subset \Gamma$,
  \end{itemize}
  such that for every $\gamma$ in
  $\Gamma$ there is some $\delta$ in $\Delta$ such that $\rho_v( \gamma s)$ is $(r,
  \varepsilon)$-proximal for every $v$ in $V$.
\end{theorem}

This result is usually stated with \emph{one} local field but the proof of the above 
extension and the following is
straightforward.

We shall also need the following Lemma

\begin{lemma}[\cite{Benoist:1997} Corollaire p.13]
  \label{lem_gammeunb}
Let  $\Gamma$ be a Zariski dense subgroup in $\mathbf{G}(k)$, $\mathbf{G}$ a
reductive group over a local field $k$. Then $\Gamma$ is unbounded if and only if  there exists an irreducible representation of $\mathbf{G}(k)$, such that $\rho( \Gamma)$ contains a proximal element.
\end{lemma}

The special case of $k=\mathbb R$ was proved in \cite{Benoist:1993}.

\subsubsection{Proximal elements and translation lengths}
\label{sec_transs}

We recall some
facts on length and translation length in $\mathbf{G}( k)$ where
$\mathbf{G}$ is a semisimple group over a local field $k$.

Let $K$ be the maximal compact subgroup of  $\mathbf{G}( k)$. This defines a norm $\Vert g \Vert_G= d_{G/K}(K,gK)$ in $\mathbf{G}( k)$ which satisfies
$$
\Vert gh\Vert_G\leqslant \Vert g\Vert_G+\Vert h\Vert_G.
$$
We also consider the translation length $\ell_G$ in $\mathbf{G}( k)$:
$$
\ell_G(g)=\inf_{h\in G}\Vert h gh^{-1}\Vert_G.
$$
Observe that the translation norm is actually independent of the choice of the maximal compact subgroup since they are all conjugated.
The translation length and norm of  $(r,\epsilon)$-proximal elements can be compared :

\begin{lemma}[compare \cite{Benoist:1997}\S~4.5]
  \label{lem_compmul}
 Let $\mathbf{G}$ be  a semisimple group over a local field $k$.  Let $\rho: \mathbf{G}(k) \rightarrow \mathrm{GL}(n, k)$ be an
  irreducible representation. Let  $\varepsilon>0$.

  Then there exist positive constants $\alpha$ and $\beta$ such that if  $\rho(g)$
  is $(r, \varepsilon)$-proximal then
  \begin{equation}
\ell_G(g)\geqslant  \alpha \Vert g\Vert_G-\beta
  \end{equation}
\end{lemma}

\proof We use classical notation and refer to \cite[p.6--8]{Benoist:1997} for precise definitions.

Any $g$ in $\mathbf{G}( k)$ is contains in a unique double coset $K
\mu(g) K$. The element $\mu(g) \in A^+ \subset A$ is called the Cartan
projection of $g$. Here we see $A^+$ as a subset of a cone $A^\times$
in some \emph{$\mathbf{R}$-vector space}.

For some integer $n$ ($n=1$ in the archimedean case) the element
$g^n$ admits a Jordan decomposition $g^n = g_e g_h g_u$ with $g_e$,
$g_h$, $g_u$ commuting, $g_e$ elliptic, $g_u$ unipotent and $g_h$
hyperbolic, \emph{i.e.} conjugated to a unique element $a \in
A^+$. We set $\lambda(g) = \frac{1}{n} a \in A^\times$.

If we fix some norm on the vector space containing the cone $A^\times$
then (up to quasi isometry constants) the norm of $\mu(g)$ is
 $\Vert g\Vert_G$ in $\mathbf{G}( k)$ and the norm of $\lambda(g)$ is $\ell_G(g)$.

Then by a result of Y. Benoist \cite{Benoist:1997}, there exists a compact subset $N_\varepsilon$ of the vector
  space containing $A^\times$ such that for every $g$ such that $\rho(g)$
  is $(r, \varepsilon)$-proximal we have
  \begin{equation*}
    \lambda(g) - \mu(g) \in N_\varepsilon.
  \end{equation*}

  The lemma follows. \qed

\vskip 0.5truecm

\subsubsection{Proof of Theorem~\ref{thm_UCC}}
\label{sec_pfthUCC}

By taking a finite index normal subgroup and projecting
(Lemmas~\ref{lem_normalsubgroup} and \ref{lem_finiteindex}) we can
make the hypothesis that $\mathbf{G}$ is the product of a semisimple
group $\mathbf{S}$ and a torus $\mathbf{T}$ and $\Gamma$ is a subgroup
of $\mathbf{S}(F) \times \mathbf{T}(F)$.

Moreover, since  length and
translation length for elements in $\mathbf{T}(F_v)$ are equal, we only need
to work with the semisimple part $\mathbf{S}$.

Finally, it suffices to  prove the existence of a finite family $F \subset
\Gamma$ and constants $A,B$ such that, for any $\gamma \in \Gamma$
\begin{equation}
  \label{eq_wanted} {\Vert} \gamma{\Vert}_S \leqslant A \sup_{f \in F} \ell_S( \gamma f) +B,
\end{equation}
where $S = \Pi_{v \in V} \mathbf{S}( F_v)$.  Indeed,  since $\Gamma$ is quasiisometrically embedded in $G$, ${\Vert} \gamma{\Vert}_{\Gamma}$
is less than ${\Vert}\gamma{\Vert}_G = {\Vert}\gamma{\Vert}_T +
{\Vert}\gamma{\Vert}_S$ and $\ell_G( \gamma) = \ell_T( \gamma) + \ell_S(
\gamma) = {\Vert}\gamma{\Vert}_T + \ell_S( \gamma)$ is less than
$\ell_\Gamma( \gamma)$  --up to quasi isometry constants-- and Inequality   \ref{eq_wanted} implies that $\Gamma$ has $U$-property.

Note that we can forget any completion $F_v$ where the subgroup
$\Gamma \subset \mathbf{S}( F_v)$ is bounded without changing the fact
that $\Gamma$ is quasiisometrically embedded. So by
Lemma~\ref{lem_gammeunb}, for each $v$ there is an irreducible
representation $\rho_v: \mathbf{S}( F_v) \rightarrow \mathrm{GL}( n_v,
F_v)$ such that $\rho_v( \Gamma)$ contains proximal elements.

Applying Theorem~\ref{coro_AMS} we find a finite family $F\subset
\Gamma$ and $r, \varepsilon$ such that for every $\gamma$ in $\Gamma$
there is some $f \in F$ such $\rho_v( \gamma f)$ is $( r,
\varepsilon)$-proximal for every $v$. Hence, as consequence of
Lemma~\ref{lem_compmul}, we have for such $\gamma$ and $f$ :
\begin{equation*}
  {\Vert} \gamma f{\Vert}_S \leqslant A \ell_S( \gamma f) +B.
\end{equation*}

This implies Inequality \eqref{eq_wanted} and concludes the proof.

\subsection{Hyperbolicity and the $U$-property}
Let $\Gamma$ be a finitely generated  group
and $S$ a set of generators. Let $d$ be its word distance and $\Vert g\Vert=d(e,g)$. We denote by
$$
\langle g,h\rangle_u=\frac{1}{2}(d(g,u)+d(h,u)-d(g,h)),
$$
the {\em Gromov product} -- based at  $u$ -- on $\Gamma$. We abbreviate $
\langle g,h\rangle_e$ by $\langle g,h\rangle$.
Observe that
\begin{eqnarray}
\langle gu,hu\rangle_u=\langle g,h\rangle_e. \label{eq:gpi}
\end{eqnarray}
Recall that $\Gamma$ is called {\em $\delta$-hyperbolic} if for all $g, h, k$ in $\Gamma$ we
have
\begin{eqnarray}
\langle g, k \rangle \geqslant  \inf  (\langle g, h \rangle ; \langle h, k \rangle) - \delta
\end{eqnarray}
and $\Gamma$ is called {\em hyperbolic} if it is $\delta$-hyperbolic for some $\delta$.  A hyperbolic group is called {\em non elementary} if it is not finite
and does not contain $\mathbb{Z}$ as a subgroup of finite index.

Then
\begin{proposition}\label{thm:hypU}
Hyperbolic groups have $U$-property.
\end{proposition}
We recall the {\em stable translation length of an element} $g$:
$$
[g]_\infty=\lim_{n\rightarrow\infty}\frac{\Vert g^n\Vert}{n}.
$$
We remark that obvioulsy
$$
[g]_\infty\leqslant \ell(g).
$$
We shall actually prove
\begin{proposition}\label{prop1}
  Let $\Gamma$ be hyperbolic. There exist a pair $u, v \in \Gamma$ and a constant
  $\alpha$
  such that for every $g$ one has
  $$\Vert g \Vert \leqslant 3\sup ( [ g ]_{\infty}, [
  g u ]_{\infty}, [ g v ]_{\infty}) + \alpha.$$
   In particular $\Gamma$ has property U.
\end{proposition}
\rmk
\begin{itemize}

 \item Let $\Gamma$ be a free group generated by some elements $u, v, w_1, \ldots
  w_n$. Then $G$ is $0$ hyperbolic. If $[ g ] \neq\Vert g \Vert$ the first
  letter of $g$ must be equal to the inverse of the last one. Multipliying
  either by $u$ or by $v$ we find a new element which is cyclically reduced :
  for this element the stable translation length and the length
  are equal. The proof of Proposition \ref{prop1} is a generalisation of
  this remark.
  \end{itemize}

\subsubsection{Almost cyclically reduced elements}
An element $g$ in $\Gamma$ is  said to be {\em almost cyclically reduced} if $\langle  g, g^{- 1} \rangle
    \leqslant  \frac{\Vert g \Vert}{3} - \delta$
We prove in this paragraph
\begin{lemma}\label{del5}
  If $g$ is almost cyclically reduced , then $$[ g ]_{\infty} \geqslant
  \frac{\Vert g \Vert}{3}.$$
\end{lemma}

The following result \cite[Lemme 1.1]{Delzant:1991} will be useful.

\begin{lemma}\label{qg1}
  Let $(x_n)$ be a finite or infinite sequence in $G$. Suppose that
  $$
  d (x_{n + 2}, x_n) \geqslant  \sup (d (x_{n + 2}, x_{n + 1}), d (x_{n
  + 1}, x_n)) + a + 2 \delta,$$
  or equivalently that
  $$
  \langle  x_{n + 2}, x_n \rangle _{x_{n + 1}} \leqslant  \frac{1}{2}\inf  (d (x_{n + 2},
  x_{n + 1}), d (x_{n + 1}, x_n)) - \frac{a}{2} - \delta.
  $$
  Then
  $$d (x_n, x_p) \geqslant  \vert n - p \vert a.$$
  \end{lemma}
This implies Lemma \ref{del5}.

\proof
  Let $x_n = g^n$. By left invariance and since $g$ is almost cyclically reduced,
  $$
  \langle  x_{n + 2}, x_n \rangle _{x_{n + 1}} = \langle  g,
  g^{- 1} \rangle  \leqslant  \frac{\Vert g \Vert}{2} - \frac{a}{2} - \delta,$$
  for $a = \frac{\Vert  g
   \Vert }{3}$. By Lemma \ref{qg1},
  $$\Vert  g^n \Vert    \geqslant    n \frac{\Vert g \Vert}{3}.$$
  The result follows
\qed
\subsubsection{Ping pong pairs}\label{def:pp}

  A {\em ping pong pair} in $ \Gamma $ is a pair of elements $u, v$ such that :
  \begin{enumerate}
    \item $\inf  ( \Vert  u \Vert  , \Vert  v  \Vert )   \geqslant    100 \delta .$
    \item $\langle  u^{\pm 1}, v^{\pm 1} \rangle \leqslant  \frac{1}{2} \inf  ( \Vert  u \Vert  , \Vert  v  \Vert )
    - 20 \delta$
    \item $\langle  u, u^{- 1} \rangle    \leqslant  \frac{\Vert  u  \Vert }{2} - 20 \delta$ and $\langle  v,
    v^{- 1} \rangle    \leqslant  \frac{\Vert  v  \Vert }{2} - 20 \delta$
  \end{enumerate}

\noindent\rmks
\begin{itemize}
\item
  A ping pong pair generates a free subgroup. This is an
  observation from \cite{Delzant:1991}. To prove this, consider a reduced
  word $w$ on the letter $u, v, u^{- 1}, v^{- 1}$.  If $x_n$ is the prefix of
  length $n$ of $w,$ the sequence $x_n$ satisfies the hypothesis of Lemma \ref{qg1}.
  \item In the present proof, the third property will not be used.
\end{itemize}
We shall prove
\begin{lemma}\label{del8}
  If $\Gamma$ is hyperbolic non-elementary, there exists a ping pong pair.
\end{lemma}

\proof  In \cite{Koubi:1998} explicit ping pong pairs are
constructed. Here is a construction whose idea goes back to F. Klein. Let
$f$ be some hyperbolic element (an element of infinite order). Replacing $f$
by a conjugate of some power, we may assume that $$
\Vert  f \Vert  = [ f ] >
1000 \delta.$$
As $\Gamma$ is not elementary, there exists a generator $a$ of $\Gamma$
which do not fix the pair of fixed points $f^+, f^-$ of $f$ on the boundary $\partial
\Gamma$ : otherwise, since the action of $\Gamma$ is topologically transitive,
$\partial \Gamma$ would be reduced to
these two points and $\Gamma$ would be elementary. Now, let us prove that for some
integer $N$, \ $(f, a f^N a^{- 1})$ is a ping pong pair. We have
$$
\lim_{N
\rightarrow + \infty} f^N = f^+\not=af^+ =\lim_{N
\rightarrow + \infty}a f^N a^{- 1}.$$
It follows that the  Gromov product $\langle f^N, a f^N a^{-
1} \rangle$ remains bounded, by the very definition of the boundary.  Hence, for  $N$ large enough, we have
$$
(f, a f^N a^{- 1})\leqslant\frac{1}{2} \inf ( \Vert  f^N \Vert ,\Vert  a f^N a^{- 1} \Vert )
- 20 \delta. $$ A similar argument  also  yields that $\langle f^N, a f^{-N} a^{-
1} \rangle$ remains bounded. Therefore $(f^N, a f^N a^{- 1})$ is a ping
pong pair for $N \gg 1$. \qed
\subsubsection{Proof of Proposition \ref{prop1}}
We first reduce this proof to the following Lemma
\begin{lemma}\label{del9}
  Let $(u, v)$ be a ping pong pair. Let $g \in  \Gamma $ such that $$
  \Vert g \Vert \geqslant 3 \sup(
  \Vert  u \Vert  , \Vert  v  \Vert ) + 100 \delta.$$
  Then one of the three elements $g, g u, g v$
  is almost cyclically reduced.
\end{lemma}

We oberve at once that Proposition \ref{prop1}  follows from Lemma \ref{del5}, \ref{del8} and \ref{del9} : choose a ping pong
  pair $u, v$ and take
  $$\alpha = 3 \sup ( \Vert  u \Vert  , \Vert  v  \Vert ) + 100 \delta.$$

\proof
Assume $g$ is not almost cyclically reduced. Then
\begin{equation}
\langle  g, g^{- 1} \rangle  \geqslant  \frac{\Vert g \Vert}{3} - \delta \geqslant   \sup( \Vert  u \Vert  , \Vert  v  \Vert )
+ 30 \delta.\label{eqdel:1}
\end{equation}
Moreover one of the following pair of inequalities hold
\begin{equation}
\langle  g^{- 1}, u^{\pm 1} \rangle  \leqslant  \frac{\Vert  u  \Vert }{2} - 10 \delta,\label{eqdel:2}
\end{equation}
or
\begin{equation}
\langle  g^{- 1}, v^{\pm 1} \rangle  \leqslant \frac{\Vert  v  \Vert }{2} - 10 \delta,\label{eqdel:3}
\end{equation}
Otherwise, by the definition of hyperbolicity we would have for some
$\varepsilon, \varepsilon' \in \left\{ \pm 1 \right\},$
$$
\langle  u^{\varepsilon}, v^{\varepsilon'} \rangle \geqslant  \frac{1}{2} \inf ( \Vert  u \Vert  ,  \Vert
v  \Vert ) - 10 \delta - \delta,$$
contradicting the  second property of the definition of ping pong pairs.

So we may assume that Inequality \ref{eqdel:2} holds. We will show that $g
u$ is almost cyclically reduced. Let $k=gu$.
By the triangle inequality $\langle  u^{- 1}, g^{} \rangle    \leqslant  \Vert  u \Vert $. Thus from
Inequality \ref{eqdel:1}, we deduce that
$$\inf  (\langle  g^{}, u^{- 1} \rangle , \langle  g, g^{- 1} \rangle ) = \langle  g, u^{-
1} \rangle.$$
Then, by the definition of hyperbolicity, we get:
\begin{equation}
\langle  g^{}, u^{- 1} \rangle    \leqslant  \langle  u^{- 1}, g^{- 1} \rangle  + \delta.\label{eqdel:3+1}
\end{equation}
Using Inequality \ref{eqdel:2} now, we have:
\begin{equation}
 \langle  g, u^{- 1} \rangle    \leqslant  \frac{\Vert  u  \Vert }{2} - 9 \delta.\label{eqdel:4}
\end{equation}
Note that $$
\langle  g, k \rangle  = \frac{1}{2} (\Vert g \Vert+\Vert  g u  \Vert -\Vert  u  \Vert )
  \geqslant    \frac{1}{2} (2 \Vert g \Vert- 2 \Vert  u  \Vert ) \geqslant  2   \Vert  u \Vert  + 100 \delta.$$
 By the triangle inequality  again
$$
\langle  k, u^{- 1} \rangle    \leqslant  \Vert  u \Vert. $$
Therefore
$$
\inf  (\langle  g, k \rangle , \langle  k, u^{- 1} \rangle ) = \langle  k, u^{- 1} \rangle.
$$
From hyperbolicity and Inequality \ref{eqdel:4}, we have
\begin{equation}
\langle  k, u^{- 1} \rangle    \leqslant  \langle  g, u^{- 1} \rangle  + \delta   \leqslant  \frac{\Vert  u
 \Vert }{2} - 8 \delta \label{eqdel:5}
\end{equation}
Applying successively Inequalities \ref{eqdel:2} and \ref{eqdel:5}, we get that
\begin{equation}
\langle  k^{- 1}, u^{- 1} \rangle  = \Vert  u \Vert  - \langle  u, g^{- 1} \rangle    \geqslant    \frac{\Vert  u  \Vert }{2} +
10 \delta   \geqslant    \langle  k, u^{- 1} \rangle  + 18 \delta.\label{eqdel:6}
\end{equation}
By hyperbolicity,
$$\inf (\langle  k, k^{- 1} \rangle , \langle  k^{- 1}, u^{- 1} \rangle )   \leqslant  \langle  k,
u^{- 1} \rangle  + \delta.
$$ Therefore, Inequalities \ref{eqdel:5} and \ref{eqdel:6} imply that
\begin{equation}
\langle  k, k^{- 1} \rangle    \leqslant  \langle  k, u^{- 1} \rangle  + \delta \leqslant  \frac{\Vert  u  \Vert }{2}
- 7 \delta\label{eqdel:7}
\end{equation}
Since  $\Vert  k \Vert    \geqslant    \Vert g \Vert - \Vert  u \Vert \geqslant  3 \Vert  u \Vert  - \Vert  u \Vert  + 100 \delta$, we finally obtain that
$k$ is almost cyclically reduced. \qed

\section{Non well displacing actions whose orbit maps are quasi isometric embeddings}
We prove in particular
\begin{proposition}
The action of ${\rm SL}(n,\mathbb Z)$ on $X_n = {\rm SL}(n,\mathbb
R)/SO(n,\mathbb R)$ is not well displacing, although, for $n \geqslant 3$, the orbit maps are quasi isometric embeddings.
\end{proposition}

The second part of this statement is a theorem of Lubotzky, Mozes and
Ragunathan \cite{Lubotzky:2000}. Note that for the action of ${\rm SL}(2, \mathbb
Z)$ on the hyperbolic plane $\mathbb H_2 = X_2$ the orbit maps are not quasi
isometries so it is obviously not well displacing since ${\rm SL}(2, \mathbb
Z)$ is a hyperbolic group (see Corollary \ref{corohyper}).

\subsection{Infinite contortion}
\label{def:contort}
We say a group $\Gamma$ has {\em infinite contortion}, if the set of conjugacy classes of powers of every non torsion element is infinite. In other words, for every non torsion element $\gamma$,   for every finite family $g_1,\ldots,g_q$ of conjugacy classes of elements of $\Gamma$, there exists $k>0$ such that
$$
\forall i\in\{1,\ldots,q\}, \gamma^k\not\in g_i.
$$

We prove
\begin{lemma}{\label{linear}}
Every residually finite group  has infinite contortion.
\end{lemma}

\proof  Let $\gamma$ be a non torsion element. Let $g_1,\ldots,g_n$ be finitely many conjugacy classes. We want to prove that there exists $k>0$ such that $\gamma^k$ belongs to no $g_i$. Since $\gamma$ is not a torsion element, we can assume that all the $g_i$ are non trivial. Let $h_i\in g_i$. Since all $h_i$ are non trivial, by residual finiteness  there exist a  homomorphism  $\phi$ in a finite groups $H$, such that  $$
\forall i, \varphi(h_i)\not=1.$$
Let $k=\Vert H\Vert$, hence $\varphi(\gamma^k)=1$. This implies that $\gamma^{k}\not\in g_i$ \qed

\subsection{Displacement function and infinite contortion }
We will prove
\begin{lemma}
Assume $\Gamma$ has infinite contortion.  Assume $\Gamma$ acts  cocompactly and properly discontinuoulsy by isometry on a space $X$. Assume furthermore that every closed bounded set in $X$ is compact. Then, for every non torsion element $\gamma$ in $\Gamma$, we have
$$
\limsup_{p\rightarrow\infty}d_X(\gamma^p)=\infty.
$$
\end{lemma}
\rmks
\begin{itemize}
\item
We should notice that the conclusion immediately fails if $\Gamma$ does not have infinite contortion. Indeed there exists an element $\gamma$ such that it powers describe only finitely many conjugacy classes of elements $g_{1}\ldots g_{q}$, and hence
$$
\limsup_{p\rightarrow\infty}d_X(\gamma^p) \le \sup_{i\in\{1,\ldots,q\}}d_{X}(g_{i})<\infty.
$$
\item It is also interesting to notice that there are groups with infinite contortion which possess elements $\gamma$ such that
$$
\liminf_{p\rightarrow\infty}d_X(\gamma^p)<\infty.
$$
Indeed, there are finitely generated linear groups  which contain elements  $\gamma$ which are conjugated to infinitely many of its powers. Hence, for such $\gamma$ we have
$$
\liminf_{p\rightarrow\infty}d_X(\gamma^p)\leqslant d_X(\gamma).
$$
Here is a simple example. We take  $\Gamma={\rm SL}(2,\mathbb Z[\frac{1}{p}])$ and
$$
\gamma=\left(
\begin{array}{cc}
1&1\\ 0&1\end{array}
\right).
$$
Then for all $n$, $\gamma^{p^n}$ is conjugated to $\gamma$.
\item 
However,
 in Paragraph \ref{bdr}, we shall give a condition - {\em bounded depth roots} (satisfied, for example,  by any group commensurable to a subgroup of ${\rm SL}(n,\mathbb Z)$) so that together with the hypothesis of previous Lemma
$$
\lim_{p\rightarrow\infty}d_X(\gamma^p)=\infty.
$$\end{itemize}
\proof  We want to  prove that
$$
\limsup_{p\rightarrow\infty}\inf_{x\in X}d(x,\gamma^p x)=\infty.
$$
Assume the contrary, then there exists
\begin{itemize}
\item a constant $R$,
\item a sequence of points $x_i$ of points in $X$,
\end{itemize}
such that for every $p$, $$d(x_i,\gamma^{p}x_i)\leqslant R.$$ Let now $K$ be a compact in $X$ such that $\Gamma.K=X$.
Let $f_i\in\Gamma$ such that $y_i=f_i^{-1}(x_i)\in K$. Then
$$d(y_i,f_i^{-1}\gamma^{p}f_i (y_i))\leqslant R.$$
Let $K_R=\{z\in X, d (z,K)\leqslant R\}$. It follows that
$$
\forall p,\ \  f_i^{-1}\gamma^{p}f_i \big( K_R\big)\cap K_R\not=0.
$$
Observe that $K_R$ is compact. By the properness of the action of $\Gamma$, we conclude that the family $\{g_i^{-1}\gamma_i^{p}g_i\}$ is finite. Hence the family of conjugacy classes of the sequence $\gamma^p$ is finite. But this contradicts infinite contortion for $\Gamma$. \qed

\begin{corollary}\label{cayl}
Assume $\Gamma$ has infinite contortion.  Let $C$ be its Cayley graph, then for $\gamma$ a non torsion element
$$
\limsup_{p\rightarrow\infty}\ell(\gamma^p)=\infty.
$$
\end{corollary}

\begin{corollary}\label{cayl2}
Assume $\Gamma$ has infinite contortion.  Let $\rho$ be a representation of
dimension $n$. Assume $\rho(\Gamma)$ contains a non trivial unipotent, then
$\rho$ is not well displacing on $X_n = {\rm SL}(n,\mathbb
R)/SO(n,\mathbb R)$.
\end{corollary}
\proof Assume $\rho$ is well displacing. Let $\gamma$ such that $\rho(\gamma)$ is a non trivial unipotent. Then for all $p$, $d_{X_n}(\gamma^p)=0$. However $\gamma$ is not a torsion element. We obtain the contradiction using the previous lemma. \qed

\subsection{Non uniform lattices}

\begin{lemma}
For $n \geqslant 3$, the action of ${\rm SL}(n,\mathbb Z)$  on $X_{n}$  is such that the orbit maps are quasi isometric embeddings. But it is not well displacing.
\end{lemma}
\proof  The group ${\rm SL}(n,\mathbb Z)$ is residually finite. Hence it has infinite contortion by Lemma \ref{linear}. The standard representation $\rho$ contains a non trivial unipotent hence it is not well displacing by Corollary \ref{cayl2}.

By a theorem of Lubotzky, Mozes and Ragunathan \cite{Lubotzky:2000} irreducible higher rank lattices $\Lambda$ are quasi-isometrically embedded in the symmetric space.
\qed

\subsection{Bounded depth roots}\label{bdr}
This section is  complementary.
We say that group $\Gamma$ has {\em bounded depth roots property}, if  for every $\gamma$ in $\Gamma$  which is a non torsion element, there exists some integer $p$, such that we have
$$
q\geqslant  p, \eta \in\Gamma\implies \eta^q\not= \gamma.
$$
Observe that ${\rm SL}(2,\mathbb Z [\frac{1}{p}])$ does not have bounded depth
roots. Note that this property is well behaved by taking subgroups and is a
property of commensurability (see Lemma~\ref{lemmasubcomens}).

We prove

\begin{proposition}\label{slnz}
The following groups have bounded depth roots property
\begin{itemize}
\item The group ${\rm SL}(n,\mathbb Z)$,
\item ${\rm SL}(n,\mathcal{O})$ where $\mathcal{O}$ is the
ring of integers of a number field $F$,
\item any subgroup of a group having bounded depth roots property or
any group commensurable to a group having this property,
\item in particular any arithmetic lattice in a archimedian Lie group.
\end{itemize}
\end{proposition}

\begin{lemma}{\sc [Bounded depth root]}\label{bdrl} Let $\Gamma$ be a group  with bounded depth roots property. Assume $\Gamma$ acts  cocompactly and properly discontinuously by isometry on a space $X$. Assume furthermore that every closed bounded set in $X$ is compact. Then, for every non torsion element $\gamma$ in $\Gamma$, we have
$$
\lim_{p\rightarrow\infty}d_X(\gamma^p)=\infty.
$$
\end{lemma}

For the proof see Pargaraph \ref{prooflembdr}. The Lemma and  Proposition above again
imply that the action of ${\rm SL}(n,\mathbb Z)$ on $X_n$ is not well displacing.

\subsubsection{Bounded depth roots property for  ${\rm SL}(n, \mathbb Z)$}
We prove the above Proposition.
\begin{lemma}{\label{bdrz}}
${\rm SL}(n,\mathbb Z)$ has bounded depth roots.
\end{lemma}

\proof
Let $A\in {\rm SL}(n,\mathbb Z)$. Let $B\in {\rm SL}(n,\mathbb Z)$. We assume there exists  $k$ such that
 $B^{k}=A$. Let $\{\lambda_j^A\}$ and $\{\lambda_j^B\}$   be the eigenvalues  of $A$ and $B$ respectively. Let
  $$
 K=\sup_{j} | \lambda_j^A |.
 $$
Then,
 $$
 \sup_{j} | \lambda_j^B | \leqslant K^\frac{1}{k}\leqslant K.
 $$
Hence  all   the coefficients of the characteristic polynomial of $B$ have a bound $K_1$ which only depends on $A$. Therefore, since these coefficients only take values in $\mathbb Z$, it follows  the  characteristic polynomials of $B$ belongs to  the finite  family
$$
\mathcal P=\{P(x)=x^n+\sum_{k=0}^{k=n-1}a_kx^k \ :\  a_k\in\mathbb Z, |  a_k |  \leqslant K_1 \}
$$
Since $\mathcal P$ is finite, there exists a constant $b>1$ such that for every root
$\lambda$ of a polynomial $P\in\mathcal P$,
$$
| \lambda| >1\implies | \lambda| \geqslant  b.
$$
Let $q\in N$ be such that $b^q$ is greater than $K$. It follows that if $B^q=A$, then all eigenvalues of $B$ have complex modulus 1. Therefore the same holds for $A$.

It follows from this discussion that we can reduce to the case where
$$
\forall i, j, |  \lambda^A_j | =1.
$$
We say such an element has {\em trivial hyperbolic part}. Note that necessarily also $B$ has trivial hyperbolic part.

We first prove that there exist a integer  $M$  depending only on  $n$, such that if $C\in {\rm SL}(n,\mathbb Z)$ has a  trivial hyperbolic part then $C^M$ is unipotent. The same argument as above show that the characteristic polynomials of elements with a trivial hyperbolic part belong a finite family of the form
$$
\mathcal P=\{P(x)=x^n+\sum_{k=0}^{k=n-1}a_kx^k \ :\   a_k\in\mathbb Z,  | a_k |  \leqslant K_2 \}.
$$
Where $K_2$ depends only on $n$.
Note that roots of polynomials belonging to $\mathcal P$ that are of complex modulus $1$ are
roots of unity. Thus we may take $M$ to be a common multiple of the orders of
those roots of unity and deduce that $C^M$ is unipotent if $C \in {\rm SL}(n,\mathbb Z)$ has a trivial hyperbolic part.

Returning to our setting we can replace $A$ by $A^M$ and $B$ by $B^M$ and  consider $A= B^k$ where both $A$ and $B$ are unipotents.
There is some rational matrix $g_0\in {\rm SL}(n,\mathbb Q)$ depending on $A$ such that $A_0 = g_0 A g_0^{-1}$ is in a Jordan form. We claim that $B_0 = g_0 B g_0^{-1}$ is a made of blocks which correspond to the Jordan blocks of $A_0$ and each such block of $B_0$ is an upper triangular unipotent matrix (This follows from observing that for unipotent matrices a matrix and its powers have the same invariant subspaces). Moreover note that the denominators of the entries of $B_0$ are bounded by some $L\in \mathbb N$ depending only on $g_0$ (and thus on $A$). By considering (some of the) entries just above the main diagonal it is easily seen that $k\le L$.
 \qed
\subsubsection{Commensurability}
We observe
\begin{lemma}\label{lemmasubcomens}
If $\Gamma$ is commensurable to a subgroup of  a group which has bounded depth roots. Then $\Gamma$ has bounded depth root.
\end{lemma}
\proof By definition a subgroup of a group having bounded depth roots has bounded depth roots. Let $G$ be a group and $H$ a subgroup having finite index $k$. Observe that for every element $g$ of $G$, we have $g^k\in H$. It follows that if $H$ has bounded depth roots, then $G$ has bounded depth roots. \qed
\subsubsection{Proof of Lemma \ref{bdrl}.}\label{prooflembdr}
\proof Let $K$ be a compact in $X$.  We first prove that
$$
\lim_{p\rightarrow\infty}\inf_{x\in K, \eta\in \Gamma}d(x,\eta^{-1}\gamma^p \eta x)=\infty.
$$
Assume the contrary, then there exists
\begin{itemize}
\item a constant $R$,
\item a  sequence of integers $p_i$ going to infinity,
\item a sequence of points $x_i$ in $K$,
\item a sequence of elements $\eta_i$ of $\Gamma$,
\end{itemize}
such that $d(x_i,\eta_i^{-1}\gamma^{p_i}\eta_i x_i)\leqslant R$. It follows that
$$
\forall i,(\eta_i^{-1}\gamma\eta_i)^{p_i} K_R\cap K_R\not=0.
$$
By the properness of the action of $\Gamma$, we conclude that the family $\{(\eta_i^{-1}\gamma\eta_i)^{p_i}\}$ is finite. But this contradicts the bounded depth root property.

We now choose the  compact $K$ such that $\Gamma.K=X$. It follows that
$$
\lim_{p\rightarrow\infty}\inf_{x\in X}d(x,\gamma^p x)=\lim_{p\rightarrow\infty}\inf_{x\in K, \eta\in \Gamma}d(\eta x,\gamma^p \eta x)=\infty.
$$
This is what we wanted to prove.
\qed

\end{document}